%%%%%%%%%% General Parameters %%%%%%%%%%%%%%%%%%%%%%%%%%%
%                                                                            %
%\baselineskip=17pt plus3pt minus1pt %increase interspace glue (default=12pt)
\baselineskip=17pt plus2pt minus1pt
%\magnification=\magstep1       %Magnify everything by 1.2
\hsize=126truemm                     %Horizontal width of text in the page
\vsize=180truemm                      %Vertical width of text on the page
\parindent=5truemm
\parskip=\smallskipamount
\mathsurround=1pt
\hoffset=2\baselineskip
\voffset=2truecm

%
%%%%%%%%%%%%%  DATE  %%%%%%%%%%%%%%%%%%%%%%%%%%%%%%%%%%%%%%
%
\def\today{\ifcase\month\or
  January\or February\or March\or April\or May\or June\or
  July\or August\or September\or October\or November\or December\fi
  \space\number\day, \number\year}
%%%%%%%%%% FONTS %%%%%%%%%%%%%%%%%%%%%%%%%%%%%%%%%%%%%%%%%%%%%%%%%%%%
%
   %romana negrita grande 14. titcap
   %romana negrita 12. seccion, teoremas,etc.
 at 10truept
\font\smallheadfont=cmr8

%%%%%%%%% Counters  %%%%%%%%%
%
\newcount\dispno      % Counter for displays/theorems/etc. 
\dispno=1\relax       % Initialized to 1
\newcount\refno       % Counter for references
\refno=1\relax        % Usually initialized to 1
\newcount\citations   % Have I cited anything yet?
\citations=0\relax    % Usually initialized to 0
\newcount\sectno      % Section number; initialized by \Chapter
\sectno=0\relax       % Initialized to zero.
\newbox\boxscratch    % A box to make calculate the width of a word
%                     % for some effects during \Section
%

%%%%%%%%%% Section Macros %%%%%%%%%%%%%%%%%%%%%%%%%%%%%%%%%%%%%%%%%%
%
% These macros are invoked when yoiu start a new Section.
% \Section is invoked with a name and nickname for the section,
% which advances the global section number. First first argument
% is the name of the section, the second is the nickname. Section
% is always at the current section number; it is incremented before
% use.
%
% For those wondering, I make a box, in boxscratch, which contains
% the string "Section m. "; then I declare a hanging indentation
% of length the width of that string, to start after the first line.
% That way, if the name of a section does not fit into a single
% line, the next lines will start to the right of "Section n.m. "
% which looks a lot nicer than having them start at the "S" of
% Section.
%
\def\Section#1#2{\global\advance\sectno by 1\relax%
\label{Section\noexpand~\the\sectno}{#2}%
\smallskip
\goodbreak
\setbox\boxscratch=\hbox{\bf Section \the\sectno.~}%
{\hangindent=\wd\boxscratch\hangafter=1
\noindent{\bf Section \the\sectno.~#1}\nobreak\smallskip\nobreak}}
%
%%%%%%%%%% end of proofs (filled square) %%%%%%%%%%%%%%%%%%%%%%%%%%
\def\sqr#1#2{{\vcenter{\vbox{\hrule height.#2pt
              \hbox{\vrule width.#2pt height#1pt \kern#1pt
              \vrule width.#2pt}
              \hrule height.#2pt}}}}
\def\square{$\mathchoice\sqr34\sqr34\sqr{2.1}3\sqr{1.5}3$}
\def\endproof{~~\hfill\square\par\medbreak}
\def\noproof{~~\hfill\square}
%
%%%%%%%%%% properly contained %%%%%%%%%%%%%%%%%%%%%%%%%%%%%%%%%%%%
%
\def\proc#1#2#3{{\hbox{${#3 \subseteq} \kern -#1cm _{#2 /}\hskip 0.05cm $}}}

%
%%%%%%%%%% normal %%%%%%%%%%%%%%%%%%%%%%%%%%%%%%%%%%%%%%%%%%%%%%%%%%%

%
\def\normalin{\hbox{\raise0.045cm \hbox
                   {$\underline{\triangleleft }$}\hskip0.02cm}}
%
%%%%%%%%%% Universalized accent %%%%%%%%%%%%%%%%%%%%%%%%%%%%%%%%%%%%
%
\def\'#1{\ifx#1i{\accent"13 \i}\else{\accent"13 #1}\fi}
%
%%% Prevents worries about accenting the `i'
%
%%%%%%%%%% MACROS %%%%%%%%%%%%%%%%%%%%%%%%%%%%%%%%%
%
% Semidirect product symbol
\def\semidirect{\rlap{$\times$}\kern+7.2778pt \vrule height4.96333pt
width.5pt depth0pt\relax\;}
%
% These macros require labels for the theorems, propositions, etc; they
% are the first parameter of them.
%
\def\prop#1#2{{\bf Proposition~\the\sectno.\the\dispno. }%
\label{Proposition\noexpand~\the\sectno.\the\dispno}{#1}\global\advance\dispno 
by 1{\it #2}\smallbreak}
\def\thm#1#2{{\bf Theorem~\the\sectno.\the\dispno. }%
\label{Theorem\noexpand~\the\sectno.\the\dispno}{#1}\global\advance\dispno
by 1{\it #2}\smallbreak}
\def\cor#1#2{{\bf Corollary~\the\sectno.\the\dispno. }%
\label{Corollary\noexpand~\the\sectno.\the\dispno}{#1}\global\advance\dispno by
1{\it #2}\smallbreak}
\def\defn{{\bf
Definition~\the\sectno.\the\dispno. }\global\advance\dispno by 1\relax}
\def\lemma#1#2{{\bf Lemma~\the\sectno.\the\dispno. }%
\label{Lemma\noexpand~\the\sectno.\the\dispno}{#1}\global\advance\dispno by
1{\it #2}\smallbreak}
\def\rmrk#1{{\bf Remark~\the\sectno.\the\dispno.}%
\label{Remark\noexpand~\the\sectno.\the\dispno}{#1}\global\advance\dispno
by 1\relax}
\def\proof{\noindent{\it Proof: }}
\def\numbeq#1{\the\sectno.\the\dispno\label{\the\sectno.\the\dispno}{#1}%
\global\advance\dispno by 1\relax}

\def\comm#1,#2{\left[#1{,}#2\right]}
\newdimen\boxitsep \boxitsep=0 true pt
\newdimen\boxith \boxith=.4 true pt 
\newdimen\boxitv \boxitv=.4 true pt
\gdef\boxit#1{\vbox{\hrule height\boxith
                    \hbox{\vrule width\boxitv\kern\boxitsep
                          \vbox{\kern\boxitsep#1\kern\boxitsep}%
                          \kern\boxitsep\vrule width\boxitv}
                    \hrule height\boxith}}
\def\square{\ \hbox{\vrule height7.5pt depth1.5pt width 6pt}\par}
\outer\def\square{\ifmmode\else\hfill\fi
   \setbox0=\hbox{} \wd0=6pt \ht0=7.5pt \dp0=1.5pt
   \raise-1.5pt\hbox{\boxit{\box0}\par}
}

\def\frac#1/#2{\leavevmode\kern.1em
              \raise.5ex\hbox{\the\scriptfont0 #1}\kern-.1em
              /\kern\.15em\lower.25ex\hbox{\the\scriptfont0 #2}}
%included but not equal
\def\incnoteq{\lower.1ex \hbox{\rlap{\raise 1ex
     \hbox{$\scriptscriptstyle\subset$}}{$\scriptscriptstyle\not=$}}}
%
%%%%%%%%%%% Commutative Diagrams Macros %%%%%%%%%%%%%%%%%%%%%%%%%%%%%%%%
%
% arrows having the label above and to the right

% arrows having the label to the left and below

% Properly contained pointing up
\def\propcontup{\bigcup\!\!\!\rlap{\kern+.2pt$\backslash$}\,\kern+1pt\vert}
%
%%%%%%%%%%  Cross Referencing Macros %%%%%%%%%%%%%%%%%%%%%%%%%%%%%%
%
%First, the macro that makes the labels
%
\def\label#1#2{\immediate\write\aux%
{\noexpand\def\expandafter\noexpand\csname#2\endcsname{#1}}}
%
% Then, a macro to test if a given macro is undefined
\def\ifundefined#1{\expandafter\ifx\csname#1\endcsname\relax}
%
%And then the macro that wrotes the labels
%
\def\ref#1{%
\ifundefined{#1}\message{! No ref. to #1;}%
 \else\csname #1\endcsname\fi}
%
% Also, macros to index the references
% automatic counting for the references
%
\def\refer#1{%
\the\refno\label{\the\refno}{#1}%
\global\advance\refno by 1\relax}
%
%% And macros to do citations automatically
%% When you make a citation to <cite-name>, it
%% creates a new control sequence, called x<cite-name>
%% and defines it as one. Then, \bibliog.mac will check
%% to see if it is defined; if so, it will print the
%% bibliography. You need \gdef, otherwise it is a local
%% definition, and it is lost when \cite is finished.
%
\def\cite#1{%
\expandafter\gdef\csname x#1\endcsname{1}%
\global\advance\citations by 1\relax
\ifundefined{#1}\message{! No ref. to #1;}%
\else\csname #1\endcsname\fi}
%
% And that's it!
%%%%%%%%%%%%%%%%%%%%%%%%%%%%%%%%%%%%%%%%%%%%%%%%%%%%%%%%%%%%%%%%%%%
%
\font\bb=msbm10 %scaled 1200     %blackboard bold
 at 8truept      %balckboard bolde
%
%
%%%%%%%%%% Blackboard Bold %%%%%%%%%%%%%%%%%%%%%%%%%%%%%%%%%%%%%%%%%%%%
%

\def\Z{\hbox{\bb Z}}

\def\Z{\hbox{\bb Z}}                     

\newread\aux
\immediate\openin\aux=\jobname.aux
\ifeof\aux \message{! No file \jobname.aux;}
\else \input \jobname.aux \immediate\closein\aux \fi
\newwrite\aux
\immediate\openout\aux=\jobname.aux

\font\smallheadfont=cmr8 at 8truept

\headline={\ifnum\pageno<2{\hfill}\else{\ifodd\pageno\rightheadline
\else\leftheadline\fi}\fi} 
\def\leftheadline{\quad\smallheadfont A. Magidin\hfil}
\def\rightheadline{\hfil\smallheadfont Absolutely closed
$\scriptstyle 2$-nil groups\quad}

\centerline{\bf ABSOLUTELY CLOSED NIL-$\bf 2$ GROUPS}
\smallskip
\centerline{Arturo Magidin}
\smallskip
{\parindent=20pt
\narrower\narrower
\noindent{\smallheadfont{Abstract. Using the description of
dominions in the variety of nilpotent groups of class at most two, we
give a characterization of which groups are absolutely closed in this
variety. We use the general result to derive an easier
characterization for some subclasses; e.g. an abelian group $\scriptstyle G$
is absolutely closed in ${\scriptstyle {\cal N}_2}$ if and only if
$\scriptstyle G/pG$ is cyclic for every prime number~$\scriptstyle p$.\par}}}
\bigskip
\medskip

\footnote{}{\noindent\smallheadfont Mathematics Subject Classification:
20E06, 20F18 (primary)}
\footnote{}{\noindent\smallheadfont Keywords: amalgam, special
amalgam, dominion, absolutely closed, nilpotent.}

The main result of this paper is a characterization of the absolutely
closed groups in the variety ${\cal N}_2$ (definitions are recalled in
\ref{prelims} below). We obtain this result by using the description
of dominions in the variety ${\cal N}_2$, and applying some ideas
D.~Saracino used in his classification of the strong
amalgamation bases for the same variety {\bf [\cite{saracino}]}.

In \ref{prelims} we will recall the main definitions and review the
notion of amalgam. In \ref{sufconds} we will recall the results of
Saracino related to his classification of amalgamation bases of~${\cal
N}_2$, and we will prove our main result. Finally, in
\ref{consequences} we will prove several reduction theorems, and
deduce some conditions which are sufficient for a group to be
absolutely closed in~${\cal N}_2$. We will also give easier to check
conditions for special classes of groups; for example, we will show
that a finitely generated abelian group is absolutely closed in~${\cal
N}_2$ if and only if it is~cyclic.

The contents of this paper are part of investigations that developed
out of the author's doctoral dissertation, which was conducted at the
University of California at Berkeley, under
the direction of Prof.~George M.~Bergman. It is my very great pleasure
to express my deep gratitude and indebtedness to Prof.~Bergman, for
his advice and encouragement throughout my graduate work and the
preparation of a prior version of this~paper, and for suggesting
\ref{classifforabs}.

\Section{Preliminaries}{prelims}

Recall that Isbell {\bf [\cite{isbellone}]} defines for a
variety $\cal C$ of algebras (in the sense of Universal Algebra) of a
fixed type~$\Omega$, and an algebra $A\in\cal C$ and subalgebra $B$ of
$A$, the {\it dominion of $B$ in $A$} to be the intersection of all
equalizers containing $B$. Explicitly,
$${\rm dom}_{A}^{\cal C}({B}) = \left\{a\in A\Bigm|\forall 
f,g\colon A\to C,\ {\rm if}\ f|_B=g|_B{\rm\ then\ }
f(a)=g(a)\right\}$$
where $C$ ranges over all algebras in~${\cal C}$, and $f,g$ are morphisms.

Also, Isbell calls an algebra $B\in{\cal C}$ {\it absolutely closed} (in~${\cal
C}$) if and only if
$$\forall A\in{\cal C}\;\;\hbox{with $B\subseteq A$,\ }
{\rm dom}_{A}^{\cal C}({B})=B.$$

For example, in the variety of semigroups, every group (when
considered as a semigroup using the forgetful functor) is absolutely
closed; this follows easily from the Zigzag Lemma~{\bf
[\cite{isbellone}]}. 

\rmrk{dependsoncontext} Note that the property of being ``absolutely
closed'' depends on the 
variety of context $\cal C$; it is common for an algebra to be
absolutely closed when considered a member of $\cal C$, and not
absolutely closed when considered as a member of a different variety
$\cal C'$.

In the variety ${\cal N}_2$ of nilpotent groups of class at most~2
(i{.}e{.} groups $G$ for which \hbox{$[G,G]\subseteq Z(G)$}) 
there are nontrivial dominions~{\bf [\cite{nildomsprelim}]}. The
precise description of dominions in this variety is recalled below. Given
that there are nontrivial dominions, 
an interesting problem is to characterize all groups that are
absolutely closed in~${\cal N}_2$.

For the remainder of this paper, every group will be assumed to lie in
${\cal N}_2$ unless otherwise specified, and all maps are assumed to
be group morphisms, unless otherwise noted. We will write all groups
multiplicatively. We will say that a group
is {\it absolutely closed} to mean it is absolutely closed in~${\cal
N}_2$. The identity element of
the group~$G$ will be denoted $e_G$, omitting the subscript if there
is no danger of ambiguity. For a group~$G$ and elements~$x$ and~$y$
in~$G$, the commutator of~$x$ and~$y$ is $[x{,}y]=x^{-1}y^{-1}xy$. The
commutator subgroup of a group $G$, denoted by $G'$ or $[G,G]$, is the
normal subgroup of $G$ generated by all $[x,y]$ with $x$, $y$ in~$G$.
More generally, given two
subsets $A$ and~$B$ of~$G$ (not necessarily subgroups), $[A,B]$
denotes the subgroup of~$G$ generated by all elements $[a,b]$, where
$a\in A$ and $b\in B$. The center of~$G$ will be denoted by~$Z(G)$. 
Any presentation of a group will be understood to be a presentation
in~${\cal N}_2$; that is, the identities of~${\cal N}_2$ will be
imposed on the group, as well as all the relations specified
in the presentation. We will use $Z$ to denote the infinite cyclic
group, which we also write multiplicatively.

In ${\cal N}_2$, since commutators are central, the commutator bracket
acts as a bilinear map from $G^{\rm ab}\times G^{\rm
ab}$ onto $[G,G]$.  In particular, for every $x,y,z\in G$, and $n\in\Z$,
$$[x,yz]=[x,y][x,z];\quad [xy,z]=[x,z][y,z];\quad
[x^n,y]=[x,y]^n=[x,y^n].$$ 

Also, given $A,B\in {\cal N}_2$, every
element of their coproduct $A\amalg^{{\cal N}_2} B$ has a unique
expression in the form $\alpha\beta\gamma$, where $\alpha\in A$,
$\beta\in B$, and $\gamma\in[A,B]$. A theorem of T. MacHenry {\bf
[\cite{machenry}]} states that the subgroup $[A,B]$ of~$A\amalg^{{\cal
N}_2}B$ is isomorphic to
the tensor product $A^{\rm ab}\otimes B^{\rm ab}$.

\def\amal#1,#2;#3{(#1{,}#2{;}#3)}
Recall that an ${\cal N}_2$-amalgam of two groups $A,C\in {\cal N}_2$
with core~$B$ consists of groups $A$, $B$, and~$C$, equipped with one to one
group morphisms
$$\eqalign{\Phi_A\colon &B\to A\cr
\Phi_C\colon &B\to C.\cr}$$
To simplify notation, we denote this situation by $(A,C;B)$. To say
that the amalgam $(A,C;B)$ is {\it (weakly) embeddable in ${\cal N}_2$} means
that there exists a group~$M$ in~${\cal N}_2$ and one-to-one group morphisms
$$\lambda_A\colon A\to M,\qquad \lambda_C\colon C\to
M,\qquad\lambda\colon B\to M$$ such that
$$\lambda_A\circ\Phi_A = \lambda \qquad\qquad
\lambda_C\circ\Phi_C=\lambda.$$

When we examine whether or not the amalgam $\amal
A,C;B$ is embeddable, the obvious candidate for $M$ is the coproduct
with amalgamation of $A$ and $C$ over~$B$, denoted
$A\amalg_B^{{\cal N}_2} C$. This coproduct is sometimes called the
${\cal N}_2$-free product with amalgamation.
We say that $\amal A,C;B$ is {\it weakly embeddable} (in~${\cal N}_2$)
if no two distinct elements of~$A$ are identified
with each other in the coproduct \hbox{$A\amalg_B^{{\cal N}_2}
C$}, and
similarly with two distinct elements of~$C$. Note that weak
embeddability does not preclude the possibility that an element $x$
of~$A\backslash B$ be identified with an element $y$ of~$C\backslash
B$ in~$A\amalg_B^{{\cal N}_2} C$. We say that $\amal A,C;B$ is {\it
strongly embeddable} (in~${\cal N}_2$)  if
there is also no identification between elements of~$A\backslash B$
and elements of~$C\backslash B$. By {\it special amalgam} we mean an amalgam
$\amal A,A';B$,
where there is an isomorphism $\psi$ between $A$ and $A'$ over $B$,
meaning that $\psi\circ\Phi_A=\Phi_{A'}$. In this case, we usually
write $\amal A,A;B$, with $\psi={\rm id}_A$ being understood.

Also, we recall that a group $G$ is said to be a {\it weak
amalgamation base} for~${\cal N}_2$ if every
amalgam with $G$ as a core is weakly embeddable in~${\cal N}_2$; it is
a {\it strong 
amalgamation base} (for~${\cal N}_2$) if every such amalgam is
strongly embeddable (in~${\cal N}_2$); and
it is a {\it special amalgamation base} for~${\cal N}_2$ if every
special amalgam with 
core~$G$ is strongly embeddable in~${\cal N}_2$. Note that a special
amalgam is always weakly~embeddable.

The connection between amalgams and dominions is via special
amalgams. Letting $A'$ be an
isomorphic copy of $A$, and
$M=A\amalg_B^{{\cal N}_2} A'$, we have that $${\rm dom}_{A}^{{\cal
N}_2}({B})=A\cap A'\subseteq M \eqno{(\numbeq{charactdom})}$$ 
where we have identified~$B$ with its common image in~$A$ and~$A'$.

The above discussion can be done in the much more general context of
an arbitrary variety $\cal C$ of algebras of a fixed type.
For a more complete discussion of amalgams in general and their
connection with dominions, see~{\bf [\cite{episandamalgs}]}.

\rmrk{siffwands} It is not hard to verify that a group $B$ is a strong
amalgamation base if and only if it is both a weak amalgamation base
and a special amalgamation base. For a proof we direct the reader
to~{\bf [\cite{episandamalgs}]}. We also note that for a group $H$ in
a variety $\cal C$, being a special amalgamation base for $\cal C$ is
equivalent to being absolutely closed in $\cal C$. Indeed, the
equality given in (\ref{charactdom}) shows that $H$ is absolutely
closed in $\cal C$ if and only if for every group $G$ containing $H$,
the special amalgam $\amal G,G;H$ is strongly embeddable, which holds
if and only if $H$ is a special amalgamation~base.

\Section{Absolutely closed groups}{sufconds}

In this section we recall the characterization of weak and strong
amalgamation bases in the variety~${\cal N}_2$, due to Saracino. Then
we will state the characterization of absolutely closed groups in this
variety. 

It will be helpful to recall a theorem about adjunction of roots
to ${\cal N}_2$-groups:

\thm{adjroots}{{\rm (Saracino, Theorem 2.1 in~{\bf
[\cite{saracino}]})} Let $G$ be a nilpotent group of class at most
two, let $m>0$, let ${\bf n}$ be an $m$-tuple of positive
integers, and let ${\bf g}$ be an $m$-tuple of elements of~$G$. Then
there exists a nilpotent group $H$ of class two, containing $G$,
and which contains an $n_i$-th root for $g_i$ ($1\leq i\leq m$) if
and only if for every $m\times m$ array $\{c_{ij}\}$ of integers such
that $n_ic_{ij}=n_jc_{ji}$ for all $i$ and~$j$, and for all
$y_1,\ldots,y_m\in G$, 
$$\hbox{if\quad} y_j^{n_j}\equiv
\prod_{i=1}^mg_i^{c_{ij}}\pmod{G'},\hbox{\quad then\quad }
\prod_{j=1}^{m}[y_j,g_j]=e.\eqno\noproof$$}

\rmrk{modulocomm} Note that \ref{adjroots} implies that we can always
adjoin $n_i$-th roots to a finite family of commutators (in fact, of central
elements). In particular, if $g\in G\in {\cal N}_2$, and $g\in
G^nG'$, then there is an extension of~$G$ which contains an $n$-th
root for $g$: since $g=x^nx'$, adjoin an $n$-th root for $x'$, and we
are done.

\thm{amalbase}{{\rm (Saracino, Theorem 3.3 in~{\bf [\cite{saracino}]})} Let
$G\in{\cal N}_2$. 
The following are equivalent:}
{\parindent=15pt
\it
\item{(i)}{$G$ is a weak amalgamation base for
${\cal N}_2$.}\par
\item{(ii)}{$G$ is a strong amalgamation base for ${\cal N}_2$.}\par
\item{(iii)}{$G$ satisfies $G'=Z(G)$, and
$\forall g\in G\;\forall n>0\; (g\in G^nG'$ or $\exists y\in
G$ and $\exists k\in\Z$ such that
$(y^n\equiv g^k \pmod{G'}\hbox{\ and\ }[y,g]\not=e))$.}
\item{(iv)}{$G$ satisfies $G'=Z(G)$, and for all $g\in G$ and
all~$n>0$, either $g$ has an $n$-th root modulo $G'$, or else $g$ has
no $n$-th root in any overgroup $K\in {\cal N}_2$ of~$G$.\noproof}\par}

We pause briefly to give some examples of groups that are strong
amalgamation bases in~${\cal N}_2$. 

{\bf Example \numbeq{exone}.} Both the dihedral group $D_8$ and the
quaternion group of 8 elements $Q$ are strong amalgamation bases.
It is clear that they lie in~${\cal N}_2$, and a routine calculation
shows that they both satisfy~{\it (iii).}

{\bf Example \numbeq{extwo}.} Analogously, any 
non-abelian group of
order $p^3$, with $p$ an odd prime, is a strong amalgamation base for
${\cal N}_2$.

\rmrk{notabs} On the other hand, we remark that a nontrivial abelian
group cannot be a strong amalgamation base in ${\cal N}_2$, since it
never satisfies $G'=Z(G)$.

Next, we recall the description of dominions in ${\cal N}_2$:

\lemma{dominion}{{\rm(See {\bf [\cite{nildomsprelim}]})} 
Let $G\in{\cal N}_2$, $H$ a subgroup of
$G$. Let $D$ be the subgroup of $G$ generated by all elements of $H$
and all elements $[x,y]^q$, where $x,y$ lie in~$G$, $q\geq0$, and 
$x^q,y^q\in H[G,G]$. Then $D = {\rm
dom}_{G}^{{\cal N}_2}(H).$\noproof}

\rmrk{alsomaier} \ref{dominion} also follows from B.~Maier's work on
amalgams of nilpotent groups; we direct the reader to~{\bf
[\cite{amalgtwo}]}.

We can now prove our main result:

\thm{bigone}{Let $G\in {\cal N}_2$. Then $G$ is absolutely closed
in~${\cal N}_2$ if and only if for all $x,y\in G$ and for all $n>0$,
one of the following holds:}
{\it\parindent=35pt
\item{(\numbeq{condtwo})}{There exist $a,b,c\in\Z$, $g_1,g_2\in G$ such that
$$\eqalign{g_1^n &\equiv x^a y^b \pmod{G'}\cr
g_2^n &\equiv x^by^c \pmod{G'}\cr}$$
and $[g_1,x][g_2,y]\not= e$; or}
\item{(\numbeq{condthree})}{There exist $a,b,c\in\Z$, $g_1,g_2\in G$ such that
$$\eqalign{g_1^n &\equiv x^a y^b \pmod{G'}\cr
g_2^n &\equiv x^{b+1}y^c \pmod{G'}.\cr}$$}\par}

\rmrk{notewhattwomeans} Note that (\ref{condtwo}) is simply the
statement that there is no extension of~$G$ which contains $n$-th
roots for both $x$ and~$y$.

\proof First, suppose that for all $x,y\in G$, and all $n>0$, either
(\ref{condtwo}) or~(\ref{condthree}) holds. 
Let $K$ be an overgroup of~$G$, and suppose that there exist $k_1,
k_2\in K$, $k_1', k_2'\in K'$ 
such that $k_1^nk_1', k_2^nk_2'\in G$. We want to show that
$[k_1,k_2]^n$ lies in~$G$. Let $x=k_1^nk_1'$ and $y=k_2^nk_2'$.

Note that since both $x$ and $y$
have $n$-th roots modulo the commutator in~$K$, there is an extension
of~$K$ which has $n$-th roots for both $x$ and~$y$ (as in
\ref{modulocomm} above). Therefore, (\ref{condtwo}) cannot hold
in~$G$. Hence, there exist $a,b,c\in \Z$, and $g_1,g_2\in G$ such
that $g_1^n\equiv x^a y^b$ and $g_2^n\equiv x^{b+1}y^c$ modulo $G'$.

Since $x=k_1^nk'$ and $y=k_2^nk'$, we have that $xk_1^{-n}$ and
$yk_2^{-n}$ are central in~$K$. In particular,
$$\Bigl[g_1k_1^{-a}k_2^{-b-1},xk_1^{-n}\Bigr]
\Bigl[g_2k_1^{-b}k_2^{-c},yk_2^{-n}\Bigr]=e.$$

On the other hand,
$$\eqalignno{[g_1k_1^{-a}k_2^{-b-1},xk_1^{-n}]
[g_2k_1^{-b}k_2^{-c},yk_2^{-n}] &=
[g_1,x][g_1,k_1^{-n}][k_1^{-a}k_2^{-b-1},x]\cr
&\qquad [k_1^{-a},k_1^{-n}][k_2^{-b-1},k_1^{-n}]\cr
&\qquad[g_2,y][g_2,k_2^{-n}][k_1^{-b}k_2^{-c},y]\cr
&\qquad[k_1^{-b},k_2^{-n}][k_2^{-c},k_2^{-n}]\cr
&= [g_1,x][g_2,y][g_1,k_1^{-n}][g_2,k_2^{-n}]\cr
&\qquad [k_1,x^{-a}][k_2,x^{-b-1}][k_1,y^{-b}][k_2,y^{-c}]\cr
&\qquad [k_2,k_1]^{(-n)(-b-1)}[k_1,k_2]^{(-b)(-n)}\cr
&= [g_1,x][g_2,y][g_1,k_1^{-n}][g_2,k_2^{-n}]\cr
&\qquad [k_1,x^{-a}y^{-b}][k_2,x^{-b-1}y^{-c}]\cr
&\qquad [k_1,k_2]^{nb - n(b+1)}\cr
&= [g_1,x][g_2,y][g_1,k_1]^{-n}[g_2,k_2]^{-n}\cr
&\qquad [k_1,g_1^{-n}][k_2,g_2^{-n}][k_1,k_2]^{-n}\cr
&= [g_1,x][g_2,y][k_1,k_2]^{-n}.\cr}$$

Therefore, $[g_1,x][g_2,y][k_1,k_2]^{-n}=e$, so
$[k_1,k_2]^n=[g_1,x][g_2,y]$; since $g_1$, $g_2$, $x$,
and~$y$ all lie in~$G$, it follows that $[k_1,k_2]^n\in G$, as~claimed.

Therefore, if $G$ satisfies the conditions,
then~$G$ is absolutely closed.

Conversely, suppose that $G$ does not satisfy the condition given.
Let $x_1,x_2\in G$, and $n>0$, such
that:\par
{\parindent=30pt
\item{(\numbeq{twofails})} For all $a,b,c\in \Z$, if $g_1,g_2\in
G$ are such that 
$$\eqalign{g_1^n&\equiv x_1^a x_2^b\pmod{G'}\cr
g_2^n&\equiv x_1^b x_2^c\pmod{G'}\cr}$$
then $[g_1,x_1][g_2,x_2]=e$; and
\item{(\numbeq{threefails})} For all $a,b,c\in \Z$, there do not exist
$g_1, g_2\in G$ such that
$$\eqalign{g_1^n&\equiv x_1^a x_2^b \pmod{G'}\cr
g_2^n&\equiv x_1^{b+1}x_2^c \pmod{G'}.\cr}$$\par}

Let $F=G\amalg^{{\cal N}_2} (Z\amalg^{{\cal N}_2}Z)$, and denote
the generators of the two copies of~$Z$ by $r_1$ and~$r_2$. Every
element of $Z\amalg^{{\cal N}_2}Z$ has a unique expression of the form
$r_1^a r_2^b [r_1,r_2]^c$. Let $N$ be the minimal normal subgroup
of~$F$ containing $x_1r_1^{-n}$ and~$x_2r_2^{-n}$. We will show that
$N\cap G=\{e\}$, and that for every $g\in G$, $g[r_1,r_2]^{-n}\notin
N$. This will prove that $G$ is not absolutely closed, by looking at
$F/N$, which contains $G$ as a subgroup, and where $[r_1,r_2]^n$ lies
in the dominion of~$G$ but not in~$G$. The proof is patterned after a proof
of Saracino (Theorem 2.1 in~{\bf [\cite{saracino}]}).

A general element of~$N$ may be written as
$$\prod_{j=1}^2\left(\prod_{k=1}^{s_j}(b_{jk}z_{jk})(x_j r_j^{-n})^{\varepsilon_{jk}}(b_{jk}z_{jk})^{-1} \right),\eqno(\numbeq{genelement})$$
where $b_{jk}\in G$, $s_j$ is a positive integer,
$\varepsilon_{jk}=\pm 1$, and $z_{jk} =
r_1^{a_{jk1}}r_2^{a_{jk2}}$. Since $F$ is nilpotent of class two, this
does indeed represent a general element of~$N$.

We may rewrite (\ref{genelement}) as follows:
$$\prod_{j=1}^2
\left(\prod_{k=1}^{s_j}\bigl[(b_{jk}z_{jk})^{-1},
(x_jr_j^{-n})^{-\varepsilon_{jk}}\bigr]
(x_jr_j^{-n})^{\varepsilon_{jk}}\right)$$
which, expanding the brackets bilinearly, becomes
$$
\prod_{j=1}^2\left(\prod_{k=1}^{s_j} [b_{jk},x_j]^{\varepsilon_{jk}}
[b_{jk},r_j^{-n}]^{\varepsilon_{jk}} [z_{jk},x_j]^{\varepsilon_{jk}}
[z_{jk},r_j^{-n}]^{\varepsilon_{jk}}\right) (x_jr_j^{-n})^{t_j}$$
where $t_j=\sum_{k=1}^{s_j}\varepsilon_{jk}$. 

Now suppose that this element is equal to an element of the form
$g[r_1,r_2]^{qn}$, for some $g\in G$, $q\in\Z$; if we write
the general expression in the form $\alpha\beta\gamma$, where
$\alpha\in G$, $\beta \in Z\amalg^{{\cal N}_2}Z$, and $\gamma\in
[G,Z\amalg^{{\cal N}_2}Z]$, then the $\beta$-factor is equal to
$r_1^{-nt_1} r_2^{-nt_2} z$, where $z$ is in the commutator of
$Z\amalg^{{\cal N}_2} Z$. But on the other hand, by uniqueness
$t_1=t_2=0$. Again by uniqueness, and using this fact, we
have:

$$\eqalign{
g & = \prod_{j=1}^2\left(\prod_{k=1}^{s_j}[b_{jk},
x_j]^{\varepsilon_{jk}}\right)\cr 
[r_1,r_2]^{qn} &= \prod_{j=1}^{2}\left(\prod_{k=1}^{s_j}
[z_{jk},r_j^{-n}]^{\varepsilon_{jk}}\right)\cr
e & = \prod_{j=1}^2\left(\prod_{k=1}^{s_j}[b_{jk},r_j^{-n}]^{\varepsilon_{jk}}
[z_{jk},x_j]^{\varepsilon_{jk}}\right).\cr}$$

Feeding in the value of $z_{jk}$ and rearranging, we have
$$\eqalignno{g &= \Bigl[\prod_{k=1}^{s_1}b_{1k}^{\varepsilon_{1k}},
x_1\Bigr] \Bigl[\prod_{k=1}^{s_2}b_{2k}^{\varepsilon_{2k}},x_2\Bigr]
&(\numbeq{generalg})\cr
[r_1,r_2]^{qn} &= [r_1,r_2]^{-n{
\Sigma\varepsilon_{2k}a_{2k1}} 
+ n{\Sigma\varepsilon_{1k}a_{1k2}}}
&(\numbeq{generalqn})\cr 
\noalign{\hbox{and}}
e &=
\prod_{j=1}^2\left[\Bigl(\prod_{k=1}^{s_j}b_{jk}^{\varepsilon_{jk}}\Bigr)^{-n}
x_1^{-\Sigma\varepsilon_{1k}a_{1kj}}
\,x_2^{-\Sigma\varepsilon_{2k}a_{2kj}}, 
r_j\right]&(\numbeq{generalmixed})\cr}$$

Now define $g_j\in G$ by $g_j = \prod_k b_{jk}^{\varepsilon_{jk}}$,
and define $c_{ij}$ by $c_{ij}=-\sum\varepsilon_{ik}a_{ikj}$. Then
(\ref{generalg}) becomes
$$g = [g_1,x_1][g_2,x_2],$$
equation (\ref{generalqn}) becomes
$$[r_1,r_2]^{qn} = [r_1,r_2]^{n(c_{21} - c_{12})},$$
and equation (\ref{generalmixed}) becomes
$$e = \left[g_1^{-n}x_1^{c_{11}}x_2^{c_{21}},r_1\right]
\left[g_2^{-n}x_1^{c_{12}}x_2^{c_{22}},r_2\right].$$

Since we know that $[G,Z\amalg^{{\cal N}_2}Z]$ is isomorphic to $G^{\rm
ab}\otimes (Z\oplus Z)$, this implies that
$$g_1^{-n}x_1^{c_{11}}x_2^{c_{21}},\; g_2^{-n}x_1^{c_{12}}x_2^{c_{22}}\in G'$$
that is,
$$\eqalign{ g_1^n &\equiv x_1^{c_{11}}x_2^{c_{21}} \pmod{G'}\cr
g_2^n &\equiv x_1^{c_{12}}x_2^{c_{22}}
\pmod{G'}.\cr}\eqno(\numbeq{congruences})$$ 
 
Now, suppose that $q=0$; that is, we are trying to find which elements
lie in $G\cap N$. Since $q=0$, it follows from (\ref{generalqn})
that $c_{21}-c_{12}=0$, that
is, that $c_{12}=c_{21}$. By (\ref{twofails}) and (\ref{congruences}),
$[g_1,x_1][g_2,x_2]=e$, and therefore, $g=e$. In
particular, $G\cap N=\{e\}$, as claimed.

Finally, suppose that $q=-1$. Then $c_{21}-c_{12}=-1$, so
$c_{12}=c_{21}+1$. But then (\ref{threefails}) says that
(\ref{congruences}) cannot occur, so there is no element $g\in G$ such
that $g[r_1,r_2]^{-n}\in N$. This proves the theorem.\endproof

In fact, we need only verify (\ref{condtwo}) and (\ref{condthree}) for
prime powers:

\cor{bigoneprimes}{Let $G\in {\cal N}_2$. Then $G$ is absolutely
closed if and only if for every $x,y\in G$, and every prime power
$p^a$, $G$ satisfies (\ref{condtwo}) or~(\ref{condthree}) with
$n=p^a$.}

\proof Necessity is immediate. To show that it is also sufficient,
note that for a given $n$, if for all $x,y\in G$, $G$ satisfies
(\ref{condtwo}) or~(\ref{condthree}), then it follows that whenever
$K$ is an overgroup of~$G$, and $k_1^nk_1', k_2^nk_2'\in G$, then
$[k_1,k_2]^n\in G$.

Let $K$ be an overgroup of~$G$, and suppose that for some $n>0$,
$k_1^nk_1',k_2^nk_2'$ both lie in~$G$. Let $n=p_1^{a_1}\cdots p_r^{a_r}$ be a
prime factorization of~$n$. Since $G$ satisfies~(\ref{condtwo})
or~(\ref{condthree}) for prime powers, it follows that
$$[k_1,k_2]^{n^2/p_i^{a_i}} =
[k_1^{n/p_i^{a_i}},k_2^{n/p_i^{a_i}}]^{p_i^{a_i}} \in G$$
for each~$i$. Let $a={\rm gcd}\{n^2/p_1^{a_1},\ldots,n^2/p_r^{a_r}\}$. Then
$[k_1,k_2]^a\in G$. But it is not hard to see that $a=n$, so
$[k_1,k_2]^n\in G$, as~claimed.\endproof

We also note the following result:

\lemma{ifdivthencondthree}{Let $G\in {\cal N}_2$, and let $n>0$. If
$x\in G^nG'$, then for all $y$ there exist $g_1,g_2\in G$ and
$a,b,c\in\Z$ such that
$$\eqalign{g_1^n&\equiv x^ay^b\pmod{G'}\cr
g_2^n&\equiv x^{b+1}y^c\pmod{G'}.\cr}$$
In particular, $G$, $n$, $x$, and~$y$ satisfy
(\ref{condthree}). Analogously, if $y\in G^nG'$, then for all $x$ we
have that $G$, $n$, $x$, and $y$ satisfy (\ref{condthree}).}

\proof Suppose that $x=r^nr'$, and $y\in G$. Let $a=b=c=0$, $g_1=e$,
and $g_2=r$. If, on the other hand, $y=s^ns'$, and $x\in G$, let
$a=c=0$, $b=-1$, $g_2=e$, and $g_1=s^{-1}$.\endproof

\cor{ifamalgthenabsclosed}{If $G\in {\cal N}_2$ is such that for every
$x\in G$ and every~$n>0$, either $x$ has an $n$-th root in~$G$
modulo~$G'$, or else $x$ does not have an $n$-th root in any ${\cal
N}_2$-overgroup of~$G$, then $G$ is absolutely closed.}

\proof Given $x,y\in G$, and $n>0$, if either $x$ or $y$ has an $n$-th
root modulo the commutator, then (\ref{condthree}) is
satisfied. Otherwise, no overgroup of~$G$ contains an $n$-th root for
either $x$ or~$y$, and hence no overgroup of~$G$ contains an $n$-th
root for {\sl both} $x$ and~$y$, so~$G$ satisfies
(\ref{condtwo}).\endproof

In particular, we deduce that any group that satisfies Saracino's
conditions is absolutely closed, which is in keeping with the fact
that any strong amalgamation base is necessarily also a special
amalgamation base.

\Section{Consequences and applications}{consequences}

First, we deduce some easy conditions from \ref{bigone} which are sufficient
for a group to be absolutely closed.

\cor{divthenac}{If $G$ is a divisible nilpotent group of class at
most~$2$, then $G$ is absolutely closed in~${\cal N}_2$.}

\proof If $G$ is divisible, then every element has an $n$-th root
modulo the commutator, so $G$ satisfies (\ref{condthree}) by
\ref{ifdivthencondthree}.\endproof

Note that any nontrivial divisible abelian group $G$ is absolutely
closed, even though it cannot be a strong amalgamation base, since the
commutator subgroup cannot equal the center.
Therefore, the class of absolutely closed
groups is strictly larger than the class of strong amalgamation bases
in ${\cal N}_2$.

Before proceeding, we will prove some reduction theorems regarding
absolutely closed groups.

If $\pi$ is a set of primes, we will say that a group $G$ is
$\pi$-divisible if every element of~$G$ has a $p$-th root in~$G$, for
every prime $p\in\pi$.
We will say that $G$ is $\pi'$-divisible if every element
of~$G$ has a $q$-th root in~$G$, for every prime $q\notin \pi$.

It is not hard to verify that for a nilpotent group $G$ of class~2,
being $\pi$-divisible is equivalent to asking that $G^{\rm ab}$ be
$\pi$-divisible. 

\thm{divisibleout}{
Let $\pi$ be a set of primes, and let $A,B\in {\cal N}_2$. Suppose
that $A$ is $\pi$-divisible, and $B$ is
$\pi'$-divisible. Then $G=A\oplus B$ is absolutely closed if and only
if both $A$ and~$B$ are.}

\proof It is easy to see that, in general, if $A\oplus B$ is
absolutely closed, then so are $A$ and~$B$.

For the converse, suppose that both $A$ and~$B$ are absolutely closed,
and let $K$ be an overgroup of~$A\oplus B$. Let $x,y\in K$, $x',y'\in
K'$, and $n>0$ be such that $x^nx', y^ny'\in A\oplus B$. We want to
show that $[x,y]^n\in A\oplus B$. Write $x^nx'= a_1\oplus b_1$, and
$y^ny' = a_2\oplus b_2$.

By \ref{bigoneprimes}, it suffices to consider the case when $n$ is a
prime power, say $n=p^{\alpha}$.

If $p\in \pi$, then $a_1^{-1}$ has an $n$-th root in~$A$.
That is, there exists $r\in A$ such that
$r^n=a_1^{-1}$.
Similarily, there exists $s\in A$ such that
$s^n=a_2^{-1}$.

Therefore,
$$\eqalign{
(rx)^n\equiv r^nx^n &\equiv a_1^{-1}(a_1\oplus b_1)\equiv
b_1\pmod{K'}\cr
(sy)^n\equiv s^ny^n &\equiv a_2^{-1}(a_2\oplus b_2)\equiv b_2
\pmod{K'}\cr}$$
so $[rx,sy]^n\in {\rm dom}_K^{{\cal N}_2}(B)$. Since $B$ is absolutely
closed, it follows that $[rx,sy]^n$ lies in~$B$. However,
$$\eqalignno{
[rx,sy]^n &= [r,s]^n[r,y]^n[x,s]^n[x,y]^n\cr
          &= [r,s]^n[r,y^ny'][x^nx',s][x,y]^n\cr
          &= [r,s]^n[r,a_2\oplus b_2][a_1\oplus b_1,s][x,y]^n\cr}$$

Since $r$ and~$s$ lie in~$A$, the first three terms on the right hand
side lie in $A\oplus
B$. Since $[rx,sy]^n\in B$, it follows that
$[x,y]^n\in A\oplus B$ as well.

If, on the other hand, $p\notin \pi$, then the argument
proceeds as above, taking roots of $b_1^{-1}$ and $b_2^{-1}$.\endproof

\cor{finiteexpred}{{\rm (Cf. Theorem 3.5 in~{\bf[\cite{saracino}]})}
If $A,B\in {\cal N}_2$ are of relatively prime exponents, then
$A\oplus B$ is absolutely closed if and only if both $A$ and~$B$ are.}

\proof If $A$ is of finite exponent $n$, then $A$ is
$\pi$-divisible, where $\pi$ is the set of all primes not occuring in 
the prime factorizaton of~$n$. The result now follows from
\ref{divisibleout}.\endproof

Recall that every abelian group $G$ may be written as $G=D\oplus
G_{\rm red}$, where $D$ is divisible and $G_{\rm red}$ is reduced. By
letting $\pi$ be the set of all primes, we obtain:

\cor{reducedpart}{An abelian group~$G$ is absolutely closed if and
only if its reduced part is absolutely closed.\noproof}

\cor{checkingsome}{If $G\in {\cal N}_2$ is $\pi$-divisible, then $G$
is absolutely closed if and only if for every $x,y\in G$ and every
prime power $n=q^a$, with $q\notin\pi$, $G$, $x$, $y$, and $n$ satisfy
(\ref{condtwo}) 
or~(\ref{condthree}).\noproof}

\cor{pparts}{If $G\in {\cal N}_2$ is a torsion group, then $G$ is
absolutely closed if and only if its $p$-parts are.\noproof}

Next we analyze what (\ref{condtwo}) and~(\ref{condthree}) mean for
finitely generated abelian groups. 

\thm{cyclicthenac}{If $G$ is cyclic, then $G$ is absolutely closed.}

\proof Let $G=\langle t\rangle$, and let $x=t^r$, $y=t^s$ be any two
elements. Let $n=p^{\alpha}$ be a prime power. We claim that $G$, $x$,
$y$, and $n$ satisfy~~(\ref{condthree}). To see this, it will suffice
to show that we can find $a$, $b$, and $c\in \Z$ such that
$p^{\alpha}|ar+bs$ and $p^{\alpha}|(b+1)r+cs$.

If $(p,r)=1$, set $b=-1$, $c=0$; then we want to find an~$a$ such
that \hbox{$p^{\alpha}|{ar}-s$}. But since $r$ is relatively prime to~$p$, as
$a$ ranges over $\Z$, $ar$ ranges over all congruence classes modulo
$p^{\alpha}$, so there is one which is congruent to $s$.

If $(p,s)=1$, we proceed similarily. Finally, suppose that
$r=p^{\delta}$, $s=p^{\gamma}$; we may assume that
$\delta,\gamma<\alpha$.

If $\delta\leq\gamma<\alpha$, then set $b=-1$, $c=0$, and
$a=p^{\gamma-\delta} + p^{\alpha-\delta}$.

And if $\gamma\leq\delta<\alpha$, then set $a=b=0$, and
let $c=-p^{\delta-\gamma} + p^{\alpha-\gamma}$.\endproof

In fact, if~$G$ is a finitely generated abelian group, then being
cyclic is also necessary for~$G$ to be absolutely closed. To prove
this, we start with a series of examples:

{\bf Example \numbeq{zsquared}.} $Z\oplus Z$ is not absolutely
closed. Indeed, let $F$ be the ${\cal N}_2$ group presented (in~${\cal
N}_2$) by
$$F= \Bigl\langle x,y\,\Bigm|\, [x,y]^4=e\Bigr\rangle;$$
then the subgroup of~$F$ generated by $x^2$ and~$y^2$ is abelian,
isomorphic to $Z\oplus Z$, but $[x,y]^2$ lies in the dominion of
$\langle x^2,y^2\rangle$, and not in the subgroup.

{\bf Example \numbeq{finitetwocyc}.} $Z/p^{a_1}Z\oplus Z/p^{a_2}Z$
with $p$ a prime, and $a_1,a_2\geq 1$, is not absolutely closed. This
time let $F$ be the ${\cal N}_2$ group presented by
$$F = \Bigl\langle x,y\,\Bigm|\,
x^{p^{a_1+1}}=y^{p^{a_2+1}}=[x,y]^{p^2}=e\Bigr\rangle$$
and let $G=\langle x^p,y^p\rangle$. Then $G\cong Z/p^{a_1}Z\oplus
Z/p^{a_2}Z$, but 
$$[x,y]^p\in {\rm dom}_F^{{\cal N}_2}(G) \setminus G.$$

{\bf Example \numbeq{zpluscyclic}.} $Z\oplus Z/p^a Z$ is not absolutely
closed, where $p$ is a prime and $a\geq 1$. Let $F$ be the ${\cal
N}_2$ group presented by
$$F=\Bigl\langle x,y\,\Bigm|\, y^{p^{a+1}} =
[x,y]^{p^2}=e\Bigr\rangle$$
and let $G=\langle x^p,y^p\rangle$. Then $G$ is isomorphic to
$Z\oplus Z/p^a Z$, and
$$[x,y]^p \in {\rm dom}_F^{{\cal N}_2}(G)\setminus G.$$

\thm{fgabelian}{A finitely generated abelian group is absolutely
closed in~${\cal N}_2$ if and only if it is cyclic.}

\proof Sufficiency is \ref{cyclicthenac}. For necessity, let $G$ be a
finitely generated abelian group, and write
$$G \cong Z^r\oplus Z/a_{1}Z\oplus\cdots\oplus Z/a_{s}Z$$
where each $a_i$ is a prime power.

If $r>1$, or $r=1$ and $s>0$, then $G$ has a direct summand which is
not absolutely closed by the examples above, hence $G$ is not
absolutely closed. If $s>1$ and there exist $i$ and~$j$ such that
$a_i$ and $a_j$ are not relatively prime, then $G$ also has a direct
summand which is not absolutely closed. All other cases (namely, $r=1$
and~$s=0$; or $r=0$ and all $a_i$ relatively prime) are cyclic
groups.\endproof

We can also prove an analogue of a result of Saracino. Recall the
following:

\thm{saracinopgroups}{{\rm (Saracino, Theorems 3.4 and~3.6 in~{\bf
[\cite{saracino}]})} Let $G$ be a nilpotent group of class 2 and
exponent $n$, where $n$ is the product of distinct primes, or twice
such a product. Then $G$ is a strong amalgamation base for~${\cal N}_2$ if
and only if $G'=Z(G)$.\noproof}

We obtain a similar result here:

\thm{pgroupsac}{Let $G$ be a nilpotent group of class two and
exponent~$n$, where $n$ is a product of distinct primes. Then $G$ is
absolutely closed if and only if $Z(G)/G'$ is cyclic.}

\proof By \ref{pparts}, we may assume that $G$ is a $p$-group, that
is, $n=p$ with $p$ a prime.  Denote the image of an element $x\in G$
in $G^{\rm ab}$ by $\overline{x}$.

Since $G^{\rm ab}$ is a $\Z/p\Z$ vector space, and $Z(G)/G'$ is a
subspace, there exist elements $\{z_i\}_{i\in I}$ and $\{b_j\}_{j\in
J}$ such that each $z_i$ lies in $Z(G)$, $\{\overline{z_i}\}$ is a
basis for $Z(G)/G'$, and $\{\overline{z_i},\overline{b_j}\}$
is a basis for $G^{\rm ab}$. Since $G$ is of exponent $p$, it follows
that $\langle z_i\,|\, i\in I\rangle$ is a direct summand of~$G$;
hence, if $|I|>1$, then $G$ is not absolutely closed. Thus, we may
assume that $|I|\leq 1$, which proves necessity.

To see sufficiency, note that if $K$ is an overgroup of~$G$, and $g\in
G$ has a $p$-th root in~$K$ modulo~$K'$, then $g$ is central in~$G$;
for if $g=r^pr'$ in~$K$, and $h\in G$, then
$$[g,h] = [r^pr',h] = [r^p,h] = [r,h^p] = [r,e] = e$$
since~$G$ is of exponent~$p$.

Also note that $G$ is $q$ divisible for any prime $q\not=p$, so it
suffices to check $p^a$-th roots. Let $K$ be any overgroup of~$G$, and
suppose that $r_1^{p^a}r_1', r_2^{p^a}r_2'\in G$, where $r_1,r_2\in
K$, $r_1',r_2'\in K'$. Write $g_1=r_1^{p^a}r_1'$, $g_2=r_2^{p^a}r_2'$. In
particular, $g_1$ and $g_2$ must be central in~$G$, hence they lie in
$\langle z_1\rangle G'$ (or in $G'$ if $|I|=0$). But then there exist
$x',y'\in G'$ such that $r_1^{p^a}r_1'x',r_2^{p^a}r_2'y'\in\langle
z_1\rangle$. Therefore,
$$[r_1,r_2]^{p^a} \in {\rm dom}_K^{{\cal N}_2}\bigl(\langle z_1\rangle\bigr) =
\langle z_1\rangle$$
since cyclic groups are absolutely closed. In particular,
$[r_1,r_2]^{p^a}\in G$, and so $G$ is absolutely closed.\endproof

Although we have proven an analogue of the ``square-free'' case of
\ref{saracinopgroups}, the ``twice a square-free number'' version does
not hold. A counterexample is:

{\bf Example \numbeq{counterextofour}.} A group $G\in {\cal N}_2$ of
exponent four, with $Z(G)/G'$ cyclic, which is not absolutely
closed. Let $G$ be presented by
$$G = \Bigl\langle x,y,z\,\Bigm|\, x^4=y^2=z^2
=[x,y]^2=[x,z]^2=[y,z]=e\Bigr\rangle.$$
Clearly, $G$ is of exponent four, and $G^{\rm ab}\cong Z/4 Z\oplus
Z/2Z\oplus Z/2Z$. Also, the center of~$G$ is generated,
modulo~$G'$, by $x^2$, so $Z(G)/G'$ is cyclic.

Let $F\in {\cal N}_2$ be presented by 
$$F=\Bigl\langle a,b,c \,\Bigm|\,
a^4=b^4=c^4=[a,b]^4=[a,c]^4=[b,c]^4=e \Bigr\rangle.$$
Then $\langle a,b^2,c^2\rangle\cong G$; yet
$$[b,c]^2\in {\rm dom}_F^{{\cal N}_2}(G)\setminus G$$
so $G$ is not absolutely closed.\endproof

In fact, we may generalize this example to show that $Z(G)/G'$ being
cyclic is no longer sufficient for finitely generated torsion
groups of exponent $p^n$, with $n>1$. Simply set 
$$G=\Bigl\langle x,y,z\,\Bigm|\,
x^{p^n}=y^p=z^p=[y,z]=[x,y]^p=[x,z]^p=e \Bigr\rangle$$
and
$$F=\Bigl\langle a,b,c\,\Bigm|\, a^{p^n} = b^{p^2} = c^{p^2} =
[a,b]^{p^2} = [a,c]^{p^2} = [b,c]^{p^2} =e\Bigr\rangle$$
and identify $G$ with the subgroup generated by $a$, $b^p$ and $c^p$.

Nevertheless, the condition that $Z(G)/G'$
be cyclic is necessary for finitely generated torsion groups:

\thm{acthenzcyclic}{Let $G\in {\cal N}_2$ be a finitely generated
(not necessarily abelian) torsion group. If $G$ is absolutely closed
in~${\cal N}_2$, then $Z(G)/G'$ is cyclic.}

\proof We may assume that $G$ is a $p$-group; suppose that $Z(G)/G'$
is not cyclic. We want to show that~$G$ is not absolutely closed. It
will suffice to show that~$G$ does not satisfy (\ref{condtwo})
or~(\ref{condthree}) for $n$ a power of~$p$.

Since $G$ is finitely generated, it is of exponent $p^a$ for some
$a>0$. Since $Z(G)/G'$ is not cyclic, there exist $x,y\in
Z(G)\setminus G'$ with
the property that if $x^ay^b\in G'$ for some integers $a,b\in\Z$, then
$x^a\in G'$ and $y^b\in G'$; simply write $Z(G)/G'$ as a sum of cyclic
groups, and let $x$ and $y$ be central elements which project to
generators of distinct cyclic summands.

Since $x$ and $y$ are both central, then (\ref{condtwo}) cannot hold
for them. Suppose then that (\ref{condthree}) holds, for $n=p^a$. Then
there exist elements $g_1,g_2\in G$, and integers $a,b,c\in \Z$, such
that
$$\eqalign{ g_1^{p^a}&\equiv x^ay^b \pmod{G'}\cr
g_2^{p^a}&\equiv x^{b+1}y^c \pmod{G'}.\cr}$$
However, $g_1^{p^a}=g_2^{p^a}=e$, hence by choice of $x$ and~$y$, we
have that $x^a$, $y^b$, $x^{b+1}$, and~$y^c$ all lie in~$G'$.

Since $G$ is a $p$-group, the orders of $x$ and $y$ modulo $G'$ are nontrivial
powers of~$p$. Therefore, we must have that $p|b$ (since $y^b\in G'$),
and that $p|b+1$ (since $x^{b+1}\in G'$). This is clearly impossible,
so $G$ does not satisfy~(\ref{condthree}). Therefore, $G$ is not
absolutely closed, as~claimed.\endproof

Using the ideas above, we can extend \ref{fgabelian} to an easy to
state characterization for all abelian groups. We start with a
technical lemma. Recall that if $G$ is an abelian group, we denote by
$nG$ the subgroup of all elements $x^n$ with $x\in G$. For an
arbitrary group~$G$, $nG$ denotes the subgroup generated by all such elements.

\lemma{abp}{For an abelian group $G$ and a prime number $p$, the
following are equivalent:}
{\parindent=20pt\it
\item{(i)} $G/pG$ is cyclic.
\item{(ii)} $G/p^aG$ is cyclic for some integer $a>0$.
\item{(iii)} $G/p^aG$ is cyclic for all integers
$a>0$.\par}

\proof Clearly (iii) implies (ii). Since
$p^aG$ is a subgroup of $pG$, it follows that $G/pG$ is a quotient of
$G/p^aG$, so (ii) implies (i). Finally, note
that for any integer $a>0$, $G/p^aG$ is an abelian group of exponent
$p^a$, hence is a direct sum of cyclic groups of orders $p^b$, with
$1\leq b\leq a$. Hence $G/pG$ is a direct sum of cyclic groups of order
$p$, with one direct summand for each direct summand in $G/p^aG$,
hence if $G/pG$ is cyclic, then so is $G/p^aG$ for each $a>0$; so (i)
implies~(iii).\endproof 

The following result was suggested by George~Bergman:

\thm{classifforabs}{Let $G$ be an abelian group (not necessarily
finitely generated). Then $G$ is absolutely closed in ${\cal N}_2$ if
and only if for every prime $p$, $G/pG$ is cyclic.}

\proof First, suppose that $G/pG$ is cyclic for each prime $p$. Let
$K$ be any overgroup of $G$, and let $x,y\in K$ be such that for some
prime $p$ and integer $a>0$, $x^{p^a}$ and $y^{p^a}$ both lie in
$G[K,K]$. We want to show that $[x,y]^{p^a}$ lies in~$G$.

By \ref{abp}, $G/p^aG$ is cyclic. Let $t\in G$ be such
that its image in $G/p^aG$ is a generator for $G/p^aG$. Let
$x',y'\in [K,K]$ be such that $x^{p^a}x',y^{p^a}y'\in G$. 

Therefore, there exist $g_1,g_2\in G$, and $r,s\in\Z$ such that
$x^{p^a}x'=t^rg_1^{-p^a}$ and $y^{p^a}y'=t^sg_2^{-p^a}$. In particular, the
elements $xg_1$ and $yg_2$ of~$K$ are such that their $p^a$-th powers
lie in $G[K,K]$; in fact, they lie in $\langle t\rangle[K,K]$. By
\ref{dominion}, $[xg_1,yg_2]^{p^a}$ lies in the dominion of~$\langle
t\rangle$. But by \ref{cyclicthenac}, the cyclic subgroup generated
by $t$ is absolutely closed, hence $[xg_1,yg_2]^{p^a}$ lies in
$\langle t\rangle$, and so in~$G$. However,
$$\eqalignno{
[xg_1,yg_2]^{p^a} &=
[x,y]^{p^a}[x,g_2]^{p^a}[g_1,y]^{p^a}[g_1,g_2]^{p^a}\cr
&= [x,y]^{p^a}[x^{p^a}x',g_2][g_1,y^{p^a}y'][g_1,g_2]^{p^a},\cr}$$
and since $g_1$, $g_2$, $x^{p^a}x'$, $y^{p^a}y'$, and
$[xg_1,yg_2]^{p^a}$ all lie in~$G$, it follows that $[x,y]^{p^a}$ also
lies in~$G$, as claimed. This shows that $G$
is absolutely~closed.

Conversely, suppose that there exists a prime $p$ such that $G/pG$ is
not cyclic. Therefore, $G/pG$ is a direct sum of more than one cyclic
group of order $p$. Let $x,y\in G$ be elements which project to
generators of distinct cyclic summands of $G/pG$. We will show that
$G$, $x$, $y$, and $n=p$ do not satisfy (\ref{condtwo}) nor
(\ref{condthree}).

Note that neither $x$ nor $y$ have $p$-th roots in~$G$, and that if a
product $x^ay^b$ has a $p$-th root in~$G$, then necessarily $p|a$ and
$p|b$.

Since $G$ is abelian, (\ref{condtwo}) cannot be satisfied. Suppose,
however, that $x$, $y$, and $p$ satisfy (\ref{condthree}). Therefore,
there exist $a,b,c\in\Z$, $g_1,g_2\in G$ such that
$$\eqalign{g_1^p &= x^ay^b\cr
g_2^p &= x^{b+1}y^c.\cr}$$

In particular, since $x^ay^b$ and $x^{b+1}y^c$ have $p$-th roots,
$p|b$ and $p|b+1$, which is clearly
impossible. Therefore, $G$, $x$, $y$, and~$n$ do not satisfy
(\ref{condthree}) either, so
$G$ cannot be absolutely closed.

This proves the theorem.\endproof

As in the case of \ref{acthenzcyclic}, when passing 
to a more general class of groups, we lose one of the implications:

\cor{sufcyclicforgen}{Let $G\in {\cal N}_2$ be a group (not
necessarily abelian). If $G/(pG)G'$ is cyclic for all primes $p$, then
$G$ is absolutely closed.}

\proof The argument above goes through, noting that instead of having
equalities $x^{p^a}x'=t^rg_1^{-p^a}$ and $y^{p^a}y'=t^sg_2^{-p^a}$, we
obtain congruences modulo $[K,K]$, which is enough for the argument to
hold.\endproof

Finally, we show that the converse of
\ref{sufcyclicforgen} does not hold:

{\bf Example \numbeq{counterexfinal}.} A group $G\in {\cal N}_2$ which
is absolutely closed, and for which $G/(3G)G'$ is not cyclic. Let $G$
be the ${\cal N}_2$ group presented by
$$G = \bigl\langle x,y,z\,\mid\,
x^3=y^3=z^3=[x,y]^3=[x,z]=[y,z]=e\bigr\rangle.$$
Then $G$ is of exponent $3$, and $Z(G)/G'$ is generated by $z$, hence
is cyclic. By \ref{pgroupsac}, $G$ is absolutely closed. Since
$3G=\{e\}$, $$G/(3G)G'\cong G/G'\cong (Z/3Z)^3,$$
so $G/(3G)G'$ is not cyclic. This shows that the condition in
\ref{sufcyclicforgen} is not necessary in general.\endproof

%% This file contains all papers and books I make reference to
%% in my write-ups. When I cite a paper, I call it up using
%% the macro \cite, which is defined in my macro package.
%% This creates a macro called \x<citation name>, and advances
%% the counter \citations by 1. At the beginning of this, I test
%% to see if \citations > 0; if so, I print ``References'',
%% and then proceed to list. I test each reference to see if
%% \x<citation name> is defined. If so, I print it, and it gets
%% assigned a number and label with \refer, also defined in my
%% macro package. It creates a label called \<citation name>,
%% whose value is the reference number, and it is writte in
%% \jobname.aux.
%
%  Reference Header, if necessary
\ifnum0<\citations{\par\bigbreak
\filbreak{\bf References}\par\frenchspacing}\fi
%
%  And then the citations, in alphabetical order
\ifundefined{xthreeNB}\else
\item{\bf [\refer{threeNB}]}{G{.} Baumslag, B.H.~Neumann,
H.~Neumann, and P.M.~Neumann. {On varieties generated by a
finitely generated group,} {\it Math.\ Z.} {\bf 86} (1964)
\hbox{93--122}. {MR:30\#138}}\par\filbreak\fi
\ifundefined{xbergman}\else
\item{\bf [\refer{bergman}]}{George M.~Bergman. {\it An Invitation to
General Algebra and Universal Constructions,''} {(Berkeley
Mathematics Lecture Notes 7, 1995)}.}\par\filbreak\fi
\ifundefined{xordersberg}\else
\item{\bf [\refer{ordersberg}]}{George M.~Bergman, {Ordering
coproducts of groups and semigroups,} {\it J.\ Algebra} {\bf 133} (1990)
no. 2, \hbox{313--339}. {MR:91j:06035}}\par\filbreak\fi
\ifundefined{xbirkhoff}\else
\item{\bf [\refer{birkhoff}]}{Garrett Birkhoff, {On the structure
of abstract algebras.} {\it Proc.\ Cambridge\ Philos.\ Soc.} {\bf
31} (1935), \hbox{433--454}.}\par\filbreak\fi
\ifundefined{xbrown}\else
\item{\bf [\refer{brown}]}{Kenneth S.~Brown, {\it Cohomology of
Groups, 2nd Edition,} {(GTM~87,
Springer Verlag,~1994)}. {MR:96a:20072}}\par\filbreak\fi
\ifundefined{xmetab}\else
\item{\bf [\refer{metab}]}{O.N.~Golovin, {Metabelian products of
groups,}
{\it Amer.\ Math.\ Soc.\ Translations series 2}, {\bf 2} (1956),
\hbox{117--131.} {MR:17,824b}}\par\filbreak\fi
\ifundefined{xhall}\else
\item{\bf [\refer{hall}]}{M.~Hall, {\it The Theory of Groups}
(Mac~Millan Company,~1959). {MR:21\#1996}}\par\filbreak\fi
\ifundefined{xphall}\else
\item{\bf [\refer{phall}]}{P.~Hall, {Verbal and marginal
subgroups,} {\it J.\ Reine\ Angew.\ Math.\/} {\bf 182} (1940)
\hbox{156--157.} {MR:2,125i}}\par\filbreak\fi
\ifundefined{xheineken}\else
\item{\bf [\refer{heineken}]}{H.~Heineken, {Engelsche Elemente der
L\"ange drei,} {\it Illinois J.\ Math.} {\bf 5} (1961)
\hbox{681--707.} {MR:24\#A1319}}\par\filbreak\fi
\ifundefined{xherman}\else
\item{\bf [\refer{herman}]}{Krzysztof~Herman, {Some remarks on
the twelfth problem of Hanna Neumann,} {\it Publ.\ Math.\ Debrecen}
{\bf 37} (1990)  no. 1--2, \hbox{25--31.} {MR:91f:20030}}\par\filbreak\fi
\ifundefined{xherstein}\else
\item{\bf [\refer{herstein}]}{I.~N. Herstein, {\it Topics in
Algebra,} (Blaisdell Publishing Co.,~1964).}\par\filbreak\fi
\ifundefined{xepisandamalgs}\else
\item{\bf [\refer{episandamalgs}]}{Peter M.~Higgins, {Epimorphisms
and amalgams,} {\it
Colloq.\ Math.} {\bf 56} no.~1 (1988) \hbox{1--17.}
{MR:89m:20083}}\par\filbreak\fi
\ifundefined{xhigmanpgroups}\else
\item{\bf [\refer{higmanpgroups}]}{Graham Higman, {Amalgams of
$p$-groups,} {\it J.\ Algebra} {\bf 1} (1964)
\hbox{301--305.} {MR:29\#4799}}\par\filbreak\fi
\ifundefined{xhigmanremarks}\else
\item{\bf [\refer{higmanremarks}]}{Graham Higman, {Some remarks
on varieties of groups,} {\it Quart.\ J.\ of Math.\ (Oxford) (2)} {\bf
10} (1959), \hbox{165--178.} {MR:22\#4756}}\par\filbreak\fi
\ifundefined{xhughes}\else
\item{\bf [\refer{hughes}]}{N.J.S.~Hughes, {The use of bilinear
mappings in the classification of groups of class~$2$,} {\it Proc.\
Amer.\ Math.\ Soc.\ } {\bf 2} (1951) \hbox{742--747.}
{MR:13,528e}}\par\filbreak\fi
\ifundefined{xisbelltwo}\else
\item{\bf [\refer{isbelltwo}]}{J.~M.~Howie and J.~R. Isbell, {
Epimorphisms and dominions II,} {\it J.\ Algebra {\bf
6}}(1967) \hbox{7--21.} {MR:35\#105b}}\par\filbreak\fi
\ifundefined{xisbellone}\else
\item{\bf [\refer{isbellone}]}{J. R. Isbell, {Epimorphisms and
dominions.} In {\it Proc.~of the Conference on Categorical Algebra, La
Jolla 1965,} 
(Lange and Springer, New
York~1966). MR:35\#105a (The statement of the
Zigzag Lemma for {\it rings} in this paper is incorrect. The correct
version is stated in~{\bf [\cite{isbellfour}]}.)}\par\filbreak\fi
\ifundefined{xisbellthree}\else
\item{\bf [\refer{isbellthree}]}{J. R. Isbell, {Epimorphisms and
dominions III,} {\it Amer.\ J.\ Math.\ }{\bf 90} (1968)
\hbox{1025--1030.} {MR:38\#5877}}\par\filbreak\fi
\ifundefined{xisbellfour}\else
\item{\bf [\refer{isbellfour}]}{J. R. Isbell, {Epimorphisms and
dominions IV,} {\it J.\ London Math.\ Soc.~(2),}
{\bf 1} (1969) \hbox{265--273.} {MR:41\#1774}}\par\filbreak\fi
\ifundefined{xjones}\else
\item{\bf [\refer{jones}]}{Gareth A.~Jones, {Varieties and simple
groups,} {\it J.\ Austral.\ Math.\ Soc.} {\bf 17} (1974)
\hbox{163--173.} {MR:49\#9081}}\par\filbreak\fi
\ifundefined{xjonsson}\else
\item{\bf [\refer{jonsson}]}{B.~J\'onsson, {Varieties of groups of
nilpotency three,} {\it Notices Amer.\ Math.\ Soc.} {\bf 13} (1966)
488.}\par\filbreak\fi
\ifundefined{xwreathext}\else
\item{\bf [\refer{wreathext}]}{L.~Kaloujnine and Marc Krasner,
{Produit complet des groupes de permutations et le probl\`eme
d'extension des groupes III,} {\it Acta Sci.\ Math.\ Szeged} {\bf 14}
(1951) \hbox{69--82}. {MR:14,242d}}\par\filbreak\fi
\ifundefined{xkhukhro}\else
\item{\bf [\refer{khukhro}]}{Evgenii I. Khukhro, {\it Nilpotent Groups
and their Automorphisms,} {(de Gruyter Expositions in Mathematics
{\bf 8}, New York 1993)}. {MR:94g:20046}}\par\filbreak\fi
\ifundefined{xkleimanbig}\else
\item{\bf [\refer{kleimanbig}]}{Yu.~G. Kle\u{\i}man, {On
identities in groups,} {\it Trans.\ Moscow Math.\ Soc.\ } 1983,
Issue 2, \hbox{63--110}. {MR:84e:20040}}\par\filbreak\fi
\ifundefined{xthirtynine}\else
\item{\bf [\refer{thirtynine}]}{L. G. Kov\'acs, {The thirty-nine
varieties,} {\it Math.\ Scientist} {\bf 4} (1979)
\hbox{113--128.} {MR:81m:20037}}\par\filbreak\fi
\ifundefined{xlamssix}\else
\item{\bf [\refer{lamssix}]}{T.Y. Lam and David B. Leep, {
Combinatorial structure on the automorphism group of~$S_6$,} {\it
Expo. Math.} {\bf 11} (1993) \hbox{289--308.}
{MR:94i:20006}}\par\filbreak\fi
\ifundefined{xlevione}\else
\item{\bf [\refer{levione}]}{F.~W. Levi, {Groups on which the
commutator relation 
satisfies certain algebraic conditions,} {\it J.\ Indian Math.\ Soc.\ New
Series} {\bf 6}(1942), \hbox{87--97.} {MR:4,133i}}\par\filbreak\fi
\ifundefined{xgermanlevi}\else
\item{\bf [\refer{germanlevi}]}{F.~W. Levi and B. L. van der Waerden,
{\"Uber eine 
besondere Klasse von Gruppen,} {\it Abhandl.\ Math.\ Sem.\ Univ.\ Hamburg}
{\bf 9}(1932), \hbox{154--158.}}\par\filbreak\fi
\ifundefined{xlichtman}\else
\item{\bf [\refer{lichtman}]}{A. L. Lichtman, {Necessary and
sufficient conditions for the residual nilpotence of free products of
groups,} {\it J. Pure and Applied Algebra} {\bf 12} no. 1 (1978),
\hbox{49--64.} {MR:58\#5938}}\par\filbreak\fi
\ifundefined{xmaxofan}\else
\item{\bf [\refer{maxofan}]}{Martin W. Liebeck, Cheryl E. Praeger, 
and Jan Saxl, {A classification of the maximal subgroups of the
finite alternating and symmetric groups,} {\it J.\ Algebra} {\bf
111}(1987), \hbox{365--383.} {MR:89b:20008}}\par\filbreak\fi
\ifundefined{xepisingroups}\else
\item{\bf [\refer{episingroups}]}{C.E. Linderholm, {A group
epimorphism is surjective,} {\it Amer.\ Math.\ Monthly\ }77
\hbox{176--177.}}\par\filbreak\fi
\ifundefined{xmckay}\else
\item{\bf [\refer{mckay}]}{Susan McKay, {Surjective epimorphisms
in classes
of groups,} {\it Quart.\ J.\ Math.\ Oxford (2),\/} {\bf 20} (1969),
\hbox{87--90.} {MR:39\#1558}}\par\filbreak\fi
\ifundefined{xmachenry}\else
\item{\bf [\refer{machenry}]}{T. MacHenry, {The tensor product and
the 2nd nilpotent product of groups,} {\it Math. Z.\/} {\bf 73}
(1960), \hbox{134--145.} {MR:22\#11027a}}\par\filbreak\fi
\ifundefined{xmaclane}\else
\item{\bf [\refer{maclane}]}{Saunders Mac Lane, {\it Categories for
the Working Mathematician,} {GTM~5},
(Springer Verlag 1971). {MR:50\#7275}}\par\filbreak\fi
\ifundefined{xabsclosed}\else
\item{\bf [\refer{absclosed}]}{Arturo Magidin, {Absolutely closed
nil-2 groups,} {{\it A.\ Universalis} to appear.}}\par\filbreak\fi
\ifundefined{xbilinearprelim}\else
\item{\bf [\refer{bilinearprelim}]}{Arturo Magidin, {Bilinear maps
and central extensions of abelian groups,} {\it Submitted.}}\par\filbreak\fi
\ifundefined{xprodvarprelim}\else
\item{\bf [\refer{prodvarprelim}]}{Arturo Magidin, {Dominions in decomposable
varieties of groups,} {\it Submitted.}}\par\filbreak\fi
\ifundefined{xmythesis}\else
\item{\bf [\refer{mythesis}]}{Arturo Magidin, {Dominions in
Varieties of Groups.} Ph.~D. thesis, University of
California at Berkeley (1998).}\par\filbreak\fi
\ifundefined{xnildomsprelim}\else
\item{\bf [\refer{nildomsprelim}]}{Arturo Magidin, {Dominions in
varieties of nilpotent groups,} {\it Comm.\ Alg.} to appear.}\par\filbreak\fi
\ifundefined{xsimpleprelim}\else
\item{\bf [\refer{simpleprelim}]}{Arturo Magidin, {Dominions in
varieties generated by simple groups,} {\it Submitted.}}\par\filbreak\fi
\ifundefined{xdomsmetabprelim}\else
\item{\bf [\refer{domsmetabprelim}]}{Arturo Magidin, {Dominions
in the variety of metabelian groups,}
{\it Submitted.}}\par\filbreak\fi
\ifundefined{xfgnilprelim}\else
\item{\bf [\refer{fgnilprelim}]}{Arturo Magidin, {Dominions in
finitely generated nilpotent groups,} {\it Comm.\
Alg.} to appear.}\par\filbreak\fi
\ifundefined{xnonsurjprelim}\else
\item{\bf [\refer{nonsurjprelim}]}{Arturo Magidin, {Nonsurjective
epimorphisms in decomposable va\-rie\-ties of\break groups,}
Submitted.}\par\filbreak\fi
\ifundefined{xwordsprelim}\else
\item{\bf [\refer{wordsprelim}]}{Arturo Magidin, {
Words and dominions,} {\it Submitted.}}\par\filbreak\fi
\ifundefined{xzabsp}\else
\item{\bf [\refer{zabsp}]}{Arturo Magidin, {$\Z$ is an absolutely
closed $2$-nil group,} {\it Submitted.}}\par\filbreak\fi
\ifundefined{xmagnus}\else
\item{\bf [\refer{magnus}]}{Wilhelm Magnus, Abraham Karras, and
Donald Solitar, {\it Combinatorial Group Theory,} 2nd Edition; (Dover
Publications, Inc.~1976). {MR:53\#10423}}\par\filbreak\fi
\ifundefined{xamalgtwo}\else
\item{\bf [\refer{amalgtwo}]}{Berthold J. Maier, {Amalgame
nilpotenter Gruppen
der Klasse zwei II,} {\it Publ.\ Math.\ Debrecen} {\bf 33}(1986),
\hbox{43--52.} {MR:87k:20050}}\par\filbreak\fi
\ifundefined{xnilexpp}\else
\item{\bf [\refer{nilexpp}]}{Berthold J. Maier, {On nilpotent
groups of exponent $p$,} {\it Journal of~Algebra} {\bf 127} (1989)
\hbox{279--289.} {MR:91b:20046}}\par\filbreak\fi
\ifundefined{xmaltsev}\else
\item{\bf [\refer{maltsev}]}{A. I. Maltsev, {Generalized
nilpotent algebras and their associated groups} (Russian), {\it
Mat.\ Sbornik N.S.} {\bf 25(67)} (1949) \hbox{347--366.} ({\it
Amer.\ Math.\ Soc.\ Translations Series 2} {\bf 69} 1968,
\hbox{1--21.}) {MR:11,323b}}\par\filbreak\fi
\ifundefined{xmaltsevtwo}\else
\item{\bf [\refer{maltsevtwo}]}{A. I. Maltsev, {Homomorphisms onto
finite groups} (Russian), {\it Ivanov. gosudarst. ped. Inst., u\v
cenye zap., fiz-mat. Nauk} {\bf 18} (1958)
\hbox{49--60.}}\par\filbreak\fi
\ifundefined{xmorandual}\else
\item{\bf [\refer{morandual}]}{S. Moran, {Duals of a verbal
subgroup,} {\it J.\ London Math.\ Soc.} {\bf 33} (1958)
\hbox{220--236.} {MR:20\#3909}}\par\filbreak\fi
\ifundefined{xhneumann}\else
\item{\bf [\refer{hneumann}]}{Hanna Neumann, {\it Varieties of
Groups,} {(Ergebnisse der Mathematik und ihrer Grenz\-ge\-biete,\/
New series, Vol.~37, Springer Verlag~1967)}. {MR:35\#6734}}\par\filbreak\fi
\ifundefined{xneumannwreath}\else
\item{\bf [\refer{neumannwreath}]}{Peter M.~Neumann, {On the
structure of standard wreath products of groups,} {\it Math.\
Zeitschr.\ }{\bf 84} (1964) \hbox{343--373.} {MR:32\#5719}}\par\filbreak\fi
\ifundefined{xpneumann}\else
\item{\bf [\refer{pneumann}]}{Peter M.~Neumann, {Splitting groups
and projectives
in varieties of groups,} {\it Quart.\ J.\ Math.\ Oxford} (2), {\bf
18} (1967),
\hbox{325--332.} {MR:36\#3859}}\par\filbreak\fi
\ifundefined{xoates}\else
\item{\bf [\refer{oates}]}{Sheila Oates, {Identical Relations in
Groups,} {\it J.\ London Math.\ Soc.} {\bf 38} (1963),
\hbox{71--78.} {MR:26\#5043}}\par\filbreak\fi
\ifundefined{xolsanskii}\else
\item{\bf [\refer{olsanskii}]}{A. Ju.~Ol'\v{s}anski\v{\i}, {On the
problem of a finite basis of identities in groups,} {\it
Izv.\ Akad.\ Nauk.\ SSSR} {\bf 4} (1970) no. 2
\hbox{381--389.}}\par\filbreak\fi
\ifundefined{xremak}\else
\item{\bf [\refer{remak}]}{R. Remak, {\"Uber minimale invariante
Untergruppen in der Theorie der end\-lichen Gruppen,} {\it
J.\ reine.\ angew.\ Math.} {\bf 162} (1930),
\hbox{1--16.}}\par\filbreak\fi
\ifundefined{xclassifthree}\else
\item{\bf [\refer{classifthree}]}{V. N. Remeslennikov, {Two
remarks on 3-step nilpotent groups} (Russian), {\it Algebra i Logika
Sem.} (1965) no.~2 \hbox{59--65.} {MR:31\#4838}}\par\filbreak\fi
\ifundefined{xrotman}\else
\item{\bf [\refer{rotman}]}{J.J. Rotman, {\it Introduction to the Theory of
Groups, 4th edition,} {(GTM~119,
Springer Verlag,~1994)}. {MR:95m:20001}}\par\filbreak\fi
\ifundefined{xsaracino}\else
\item{\bf [\refer{saracino}]}{D. Saracino, {Amalgamation bases for
nil-$2$ groups,} {\it Alg.\ Universalis\/} {\bf 16} (1983),
\hbox{47--62.} {MR:84i:20035}}\par\filbreak\fi
\ifundefined{xscott}\else
\item{\bf [\refer{scott}]}{W. R. Scott, {\it Group Theory,} (Prentice
Hall,~1964). {MR:29\#4785}}\par\filbreak\fi
\ifundefined{xsmelkin}\else
\item{\bf [\refer{smelkin}]}{A. L. \v{S}mel'kin, {Wreath products and
varieties of groups,} (Russian) {\it Dokl.\ Akad.\ Nauk S.S.S.R.\/} {\bf
157} (1964), \hbox{1063--1065} Transl.: {\it Soviet Math.\ Dokl.\ } {\bf
5} (1964), \hbox{1099--1011}. {MR:33\#1352}}\par\filbreak\fi
\ifundefined{xstruikone}\else
\item{\bf [\refer{struikone}]}{Ruth Rebekka Struik, {On nilpotent
products of cyclic groups,} {\it Canadian J.\ Math.}
{\bf 12} (1960)
\hbox{447--462}. {MR:22\#11028}}\par\filbreak\fi
\ifundefined{xstruiktwo}\else
\item{\bf [\refer{struiktwo}]}{Ruth Rebekka Struik, {On nilpotent
products of cyclic groups II,} {\it Canadian J.\ Math.}
{\bf 13} (1961) \hbox{557--568.}
{MR:26\#2486}}\par\filbreak\fi
\ifundefined{xvlee}\else
\item{\bf [\refer{vlee}]}{M. R. Vaughan-Lee, {Uncountably many
varieties of groups,} {\it Bull.\ London Math.\ Soc.} {\bf 2} (1970)
\hbox{280--286.} {MR:43\#2054}}\par\filbreak\fi
\ifundefined{xweibel}\else
\item{\bf [\refer{weibel}]}{Charles Weibel, {\it Introduction to
Homological Algebra,} (Cambridge University
Press~1994). {MR:95f:18001}}\par\filbreak\fi 
\ifundefined{xweigelone}\else
\item{\bf [\refer{weigelone}]}{T. S. Weigel, {Residual properties
of free groups,} {\it J.\ Algebra} {\bf 160} (1993)
\hbox{14--41.} {MR:94f:20058a}}\par\filbreak\fi
\ifundefined{xweigeltwo}\else
\item{\bf [\refer{weigeltwo}]}{T. S. Weigel, {Residual properties
of free groups II,} {\it Comm.\ in Algebra} {\bf 20}(5) (1992)
\hbox{1395--1425.} {MR:94f:20058b}}\par\filbreak\fi
\ifundefined{xweigelthree}\else 
\item{\bf [\refer{weigelthree}]}{T. S. Weigel, {Residual Properties
of free groups III,} {\it Israel J.\ Math.\ } {\bf 77} (1992)
\hbox{65--81.} {MR:94f:20058c}}\par\filbreak\fi
\ifundefined{xzstwo}\else
\item{\bf [\refer{zstwo}]}{Oscar Zariski and Pierre Samuel,
{\it Commutative Algebra, Volume
II,} (Springer-Verlag~1976). {MR:52\#10706}}\par\filbreak\fi
\ifnum0<\citations\nonfrenchspacing\fi

\bigskip

{\it
\obeylines
\noindent Arturo Magidin
\noindent Oficina 112
\noindent Instituto de Matem\'aticas
\noindent Universidad Nacional Aut\'onoma de M\'exico
\noindent 04510 Mexico City, MEXICO
\noindent e-mail: magidin@matem.unam.mx
}

\vfill
\eject
\immediate\closeout\aux
\end